\documentclass[12pt,leqno]{article}
\usepackage{amsmath,amsfonts,amssymb,amsthm}

\def\q \m#1#2{{\raise1pt\hbox{$#1$}\kern-1pt\big/
               \kern-1pt\raise-1pt\hbox{$#2$}}}

\def\R{{\mathbb R}}
\def\Z{{\mathbb Z}}

\def\C{{\mathbb C}}
\def\Aut{{\rm Aut}}
\def\wt{{\rm wt}}
\def\End{{\rm End}}
\def \pf {\noindent {\bf Proof :} \,}
\def \qed{\mbox{ $\square$}}

\def\1{{\bf 1}}
\def\o{{\omega}}
\def\l{{\lambda}}
\def\<{\langle}
\def\>{\rangle}
\def\wi{\widetilde}
\def\Res{{\rm Res}}
\def\sym{{\rm Sym}}
\def\a{\alpha}
\def\b{\beta}
\def\M{{\underline{M}}}
\def\w{{\omega}}
\def\O{{\cal O}}
\def\Hom{{\rm Hom}}

\newfam\msbfam

\font\twelmsb=msbm10 at 12pt 
\font\sevenmsb=msbm10 at 7pt \font\fivemsb=msbm10 at 5pt

\textfont\msbfam=\twelmsb
\scriptfont\msbfam=\sevenmsb \scriptscriptfont\msbfam=\fivemsb
\newtheorem{thm}{Theorem}[section]

\newtheorem{lem}[thm]{Lemma}
\newtheorem{lemma}[thm]{Lemma}
\newtheorem{prop}[thm]{Proposition}
\newtheorem{cor}[thm]{Corollary}

\theoremstyle{remark}
\newtheorem{rem}[thm]{Remark}
\theoremstyle{definition}
\newtheorem{defn}[thm]{Definition}

\newcommand{\cC}{{\cal C}}

\newcommand{\g}{{\mathfrak g}}
\newcommand{\h}{{\mathfrak h}}

\newcommand{\comment}[1]{}

\begin{document}

\renewcommand{\theequation}{\thesection.\arabic{equation}}
\setcounter{equation}{0}

\centerline{\Large {\bf $K$-theory associated to vertex operator
algebras}} \vskip 7mm \centerline{\bf Chongying
DONG\footnote{Mathematics Department, University of California,
Santa Cruz, CA 95064, U.S.A. (dong@math.ucsc.edu), Partially
supported by NSF grants and a research grant from the Committee on
Research, UC Santa Cruz.},
 Kefeng   LIU\footnote{Center of Mathematical Sciences, Zhejiang University and Department of Mathematics, UCLA,CA 90095-1555,
USA (liu@math.ucla.edu). Partially supported by an NSF grant.},
Xiaonan MA\footnote{Centre de Math\' ematiques, UMR 7640 du CNRS,
Ecole Polytechnique, 91128 Palaiseau Cedex, France
(ma@math.polytechnique.fr)}, Jian ZHOU\footnote{Department of
Mathematics, Qinghua University, Beijing, China
(jzhou@math.tsinghua.edu.cn)}} \vskip 5mm

{\bf Abstract.}  We introduce two $K$-theories, one for vector
bundles whose fibers are modules of vertex operator algebras,
another for vector bundles whose fibers are modules of associative
algebras. We verify the cohomological properties of these
$K$-theories, and construct a natural homomorphism from the VOA
K-theory to the associative algebra K-theory.

\section{ \normalsize Introduction}
\setcounter{equation}{0}

$\quad$ Since its introduction by Grothendieck, Atiyah and
Hirzebruch, $K$-theory has found many applications in algebraic
geometry, topology and differential geometry. K-theories in
different settings lead to the
Grothendieck-Hirzebruch-Riemann-Roch theorem and the Atiyah-Singer
index theory. Originally such theory was developed starting from
vector bundles. Note that one can regard vector spaces, which are
fibers of vector bundles, simply as $\mathbb C$-modules. It is
natural to consider bundles of modules over other algebras.

The theory of vertex operator algebras has been very much
developed in the last eighteen years. Playing important roles in
the study of elliptic genus and Witten genus, the highest weight
representations of Heisenberg and affine Kac-Moody algebras
provide important examples of vertex operator algebras. In this
note we introduce a $K$-theory for vector bundles of modules of
vertex operator algebras, and a K-theory for vector bundles of
modules of associative algebras.  We verify the cohomology theory
properties of these $K$-theories -- the exact sequences. We also
give  a natural homomorphism from the vertex operator algebra
K-group to the associative algebra K-group when the associative
algebra is the Zhu's algebra of the vertex operator algebra.

We follow closely the construction of topological K-theory by
Atiyah in \cite{A} to verify the cohomological properties of the
two new K-groups. We need certain deep property of vertex operator
algebras, such as the existence of the nondegenerate symmetric
bilinear form, to get exact sequences. Our original motivation for
this work is to understand elliptic cohomology which only has
homotopy-theoretical construction. We plan to study the
applications of such K-theories, in particular their relationship
to the elliptic genus and elliptic cohomology.

\section{ \normalsize Vertex operator algebras and modules}
\setcounter{equation}{0}

$\quad$ We recall the definitions of vertex operator algebras and
their modules (cf. \cite{B86}, \cite{DLM1}, \cite{FLM},
\cite{FHL},\cite{Z}) in this section. Rational vertex operator
algebras and contragredient modules are discussed. We will also
introduce the category $\O_V$ of ordinary $V$-modules for the
purpose of $K$-theory developed in this paper. The Heisenberg
vertex operator algebras and vertex operator algebras associated
to the highest weight representations for affine Kac-Moody
algebras are reviewed. They have played important roles in the
study of elliptic genus and Witten genus in the literature
(cf.\cite{W}).

\subsection{ \normalsize Vertex operator algebras and modules}
$\quad$ Let $z,$ $z_0,$ $z_1,$ $z_2$ be commuting formal
variables. We shall use the formal power series
$\delta(z)=\sum_{n\in{\Z}}z^n$ which is formally the expansion of
the $\delta$-function at $z=1.$ For properties of the
$\delta$-function see e.g. \cite{FLM}, \cite{FHL} and \cite{DL}.

A {\em vertex  operator algebra} is a ${\Z}$-graded vector space:
\begin{eqnarray} \label{0a2}
V=\bigoplus_{n\in{\Z}}V_n; \ \ \ {\rm for}\ \ \ v\in V_n,\ \
n={\rm wt}\,v;
\end{eqnarray}
such that $\dim V_n<+\infty$ for all $n\in\Z$ and $V_n=0$ if $n$
is sufficiently small; equipped with a linear map
 \begin{align}\label{0a3}
& V \to (\mbox{End}\,V)[[z,z^{-1}]] \\
& v\mapsto Y(v,z)=\sum_{n\in{\Z}}v_nz^{-n-1}\ \ \ \  (v_n\in
\mbox{End}\,V)\nonumber
\end{align}
and with two distinguished vectors ${\bf 1}\in V_0,$ $\omega\in
V_2$ satisfying the following conditions for $u, v \in V$:
\begin{align} \label{0a4}
& u_nv=0\ \ \ \ \ {\rm for}\ \  n\ \ {\rm sufficiently\ large};  \\
& Y({\bf 1},z)=1;  \\
& Y(v,z){\bf 1}\in V[[z]]\ \ \ {\rm and}\ \ \ \lim_{z\to
0}Y(v,z){\bf 1}=v;
\end{align}
and there exists a nonnegative integer $n$ depending on $u,v$ such
that
\begin{align} \label{0a7}
&(z_1-z_2)^n[Y(u,z_1), Y(v,z_2)]=0;
\end{align}
Write
\begin{eqnarray} \label{0a9}
L(n)=\omega_{n+1}\ \ \ \mbox{for}\ \ \ n\in{\Z}, \ \ \
\mbox{i.e.},\ \ \ Y(\omega,z)=\sum_{n\in{\Z}}L(n)z^{-n-2},
\end{eqnarray}
then
\begin{align}
&[L(m),L(n)]=(m-n)L(m+n)+\frac{1}{12}(m^3-m)\delta_{m+n,0}(\mbox{rank}\,V),
\end{align}
for $m, n\in {\Z},$ and
\begin{align} \label{0a10}
&L(0)v=nv=(\mbox{wt}\,v)v \ \ \ \mbox{for}\ \ \ v\in V_n\
(n\in{\Z}); \\
&\frac{d}{dz}Y(v,z)=[L(-1),Y(u,z)]=Y(L(-1)v,z).
\end{align}
This completes the definition. We denote the vertex operator
algebra just defined by $(V,Y,{\bf 1},\omega)$ (or briefly, by
$V$). The series $Y(v,z)$ are called vertex operators.

\begin{rem}\label{rem1}{\rm Any commutative associative algebra $V$ over $\C$
is a vertex operator algebra with $\1=1,$ $\omega=0$ and
$Y(u,z)v=uv$ for $u,v\in V.$ In particular, $\C$ itself is a
vertex operator algebra.}
\end{rem}

An {\em automorphism} $g$ of the vertex operator algebra $V$ is a
linear automorphism of $V$ preserving ${\bf 1}$ and $\omega$  such
that the actions of $g$ and $Y(v,z)$ on $V$ are compatible in the
sense that $gY(v,z)g^{-1}=Y(gv,z)$ for $v\in V.$ Then $gV_n\subset
V_n$ for $n\in\Z.$ The group of all automorphisms of the vertex
operator algebra $V$ is denoted by $\Aut(V)$.

We now  define admissible modules and ordinary modules for vertex
operator algebras. An {\em admissible} $V$-module
$$M=\bigoplus_{n=0}^{\infty}M(n)$$
is a $\Z$-graded vector space with the top level $M(0)\ne 0$
equipped with a linear map
\begin{align}\label{0a15}
&V\longrightarrow (\End\,M)[[z,z^{-1}]]\\
&v\longmapsto\displaystyle{ Y_M(v,z)=\sum_{n\in\Z}v_nz^{-n-1}\ \ \
(v_n\in \End\,M)} \nonumber
\end{align}
which satisfies the following conditions; for $u,v\in V,$ $w\in
M$, $n\in \Z,$
\begin{eqnarray}\label{0a16}\begin{array}{lll}
& &u_nw=0\ \ \
\mbox{for}\ \ \ n\gg 0,\\
& &Y_M({\bf 1},z)=1,\end{array}
\end{eqnarray}
\begin{equation}\label{jacm}
\begin{array}{c}
\displaystyle{z^{-1}_0\delta\left(\frac{z_1-z_2}{z_0}\right)
Y_M(u,z_1)Y_M(v,z_2)-z^{-1}_0\delta\left(\frac{z_2-z_1}{-z_0}\right)
Y_M(v,z_2)Y_M(u,z_1)}\\
\displaystyle{=z_2^{-1}\delta\left(\frac{z_1-z_0}{z_2}\right)
Y_M(Y(u,z_0)v,z_2)},
\end{array}
\end{equation}
which is called the Jacobi identity. Here all binomial expressions
$(z_i-z_j)^n$ are to be expanded in nonnegative integral powers of
second variable $z_j.$ This identity is interpreted algebraically
as follows: if this identity is applied to a single vector of $M$
then the coefficient of each monomial in $z_0,z_1,z_2$ is a finite
sum in $M;$
\begin{eqnarray}\label{0a18}
u_mM(n)\subset M(\wt u-m-1+n)
\end{eqnarray}
if $u$ is homogeneous. We denote the admissible $V$-module by
$M=(M,Y_M).$

Homomorphisms and isomorphisms of admissible modules are defined
as expected.

\begin{rem} Let $(M,Y_M)$ be an admissible $V$-module. Then
$L(-1)$-derivation property
\begin{equation}\label{moti}
Y_M(L(-1)v,z)=\frac{d}{dz}Y_M(v,z)
\end{equation}
holds. Moreover, the component operators of $Y_M(\o,z)$ generate a
copy of the Virasoro algebra of central charge ${\rm rank}\,V$
(see \cite{DLM1}).
\end{rem}

\begin{defn} Let $G$ be a subgroup of $\Aut(V)$ and $M=(M,Y_M)$ an admissible
$V$-module. We say that $G$ acts on $M$ as automorphisms if $M$ is
a $G$-module and $gY_M(v,z)g^{-1}=Y_M(gv,z)$ on $M$ for all $g\in
G, v\in V.$
\end{defn}

\begin{rem} In general $G$ does not act on an admissible module. Assume
that $V=\sum_{n=0}^{\infty}V_n$ with $V_0=\C\1.$ Then $V_1$ is a
Lie algebra. Let $N$ be the subgroup of $\Aut(V)$ generated by
$e^{u_0}$ for $u\in V_1.$ The group  $N$ acts on any admissible
module.
\end{rem}

An ({\em ordinary})  $V$-module is an admissible $V$-module $M$
which carries a $\C$-grading induced by the spectrum of $L(0).$
That is, we have
$$M=\bigoplus_{\lambda \in{\C}}M_{\lambda} $$
where $M_{\l}=\{w\in M|L(0)w=\l w\}.$ Moreover we require that
$\dim M_{\l}$ is finite and for fixed $\l,$ $M_{n+\l}=0$ for all
small enough integers $n.$ Let $P(M)=\{\l\in\C |M_{\l}\ne 0\}.$ An
element $\l\in P(M)$ is called a weight of $M.$

It is easy to show that an ordinary module is an admissible
module.

\begin{defn} A vertex operator algebra $V$ is called {\em rational} if any
admissible $V$-module is a direct sum of irreducible admissible
$V$-modules.
\end{defn}

It is  proved in \cite[Theorem 8.1]{DLM2} (also see \cite{Z}) that
if $V$ is a rational vertex operator algebra then every
irreducible admissible $V$-module is an ordinary $V$-module and
$V$ has only finitely many irreducible admissible  modules up to
isomorphism.

\begin{defn} For each $\l\in\C$ we set
$D(\l)=\{\l+n|0\leq n\in\Z\}.$ We define a category $\O_V$ of
ordinary $V$-modules such that for each $M\in \O_V$ there are
finitely many complex numbers $\l_1,...,\l_s$ such that
$P(M)\subset \cup_{i=1}^s D(\l_i).$
\end{defn}

Note that if $M$ is an irreducible $V$-module then there exists
$\lambda\in\C$ such that $M=\bigoplus_{n\geq 0}M_{\l+n}.$ Thus any
irreducible $V$-module lies in the category $\O_V.$

The definition of category $\O_V$ is very similar to the category
$\O$ in the representation theory for Kac-Moody algebras (cf.
\cite{K}). The purpose for such definition will be clear from the
discussion in the next section.

\begin{rem}\label{raa}{\rm (1) If $V$ is rational then the category $\O_V$ is exactly
the category of ordinary $V$-modules.

(2) If $U$ and $V$ are two vertex operator algebras, there is a
functor from $\O_U\times \O_V \to \O_{U\otimes V}$ such that
$(M,N)$ is mapped to $M\otimes N.$ }
\end{rem}

Next we recall the notion of contragredient module from
\cite{FHL}. Let $M=\bigoplus_{\lambda \in{\C}}M_{\lambda} $ be a
$V$-module and $M'=\bigoplus_{\lambda \in{\C}}M_{\lambda}^*$ the
graded dual. We denote the natural pairing on $M'\times M$ by
$(w',w)$ for $w'\in M'$ and $w\in M.$ Then $(M',Y)$ is also a
$V$-module such that
\begin{eqnarray} \label{0a14}
(Y(u,z)w',w)=(w', Y(e^{zL(1)}(-z^{-2})^{L(0)}v,z^{-1})w)
\end{eqnarray}
for $u\in V,$ $w'\in M'$ and $w\in W$ (see \cite[Theorem
5.2.1]{FHL}). Moreover, $M'$ is irreducible if and only if $M$ is
irreducible.

Let $W$ be a $V$-module. A bilinear form on $W$ is called {\em
invariant} if
$$(Y(u,z)w_1,w_2)=(w_1, Y(e^{zL(1)}(-z^{-2})^{L(0)}v,z^{-1})w_2)$$
for $w_i\in W$ and $v\in V.$

\begin{lem}\label{l2.1} For any $V$-module $M,$ $M\oplus M'$ has a natural
nondegenerate symmetric invariant bilinear form defined by
$$(u+u',w+w')=(u,w')+(w,u')$$
for any $u,w\in M$ and $u',w'\in M'.$ In particular, any
$V$-module can be embedded into a module with a nondegenerate
symmetric invariant bilinear form.
\end{lem}

\subsection{ \normalsize Examples}

$\quad$ In order to discuss examples of vertex operator algebra
bundles in the next few sections we recall some well-known vertex
operator algebras.

(a) Heisenberg vertex operator algebra $M(1)$ of dimension $d$
(cf. \cite{FLM}). Let ${\h}$ be complex vector space of dimension
$d$ with a nondegenerate symmetric bilinear from $(,).$ Viewing
${\h}$ as an abelian Lie algebra and consider the corresponding
affine Lie algebra $\hat{\h}={\h}\otimes \C[t,t^{-1}]\oplus \C K$
with bracket
$$[x\otimes t^m,y\otimes t^n]=\delta_{m+n,0}(x,y)K, $$
$$[K,\hat{\mathfrak h}]=0,$$
where $x(m)=x\otimes t^m$ for $x\in \mathfrak h$ and $m\in \Z.$
Form the induced module
$$M(1)=U(\hat{\h})\otimes_{U(\h\otimes \C[t]\oplus \C K)}\C,$$
where $\h\otimes \C[t]$ acts trivially on $\C$ and $K$ acts
as 1. Let $\{\alpha_1,...,\alpha_d\}$ be an orthonormal basis of
${\h}.$ Then $M(1)$ is isomorphic linearly to the symmetric
algebra
$$S(\h\otimes t^{-1}\C[t^{-1}])=\C[\alpha_i(-n)|i=1,...,d,n>0].$$
Set $\1=1$ and $\omega=\frac{1}{2}\sum_{i=1}^d\alpha_i(-1)^2.$
Then $M(1)$ is the Heisenberg vertex operator algebra with vacuum
$\1,$ Virasoro element $\omega$ (cf. \cite[Chapter 8]{FLM}). Let
$O(\h)$ be the orthogonal group of $\h$. The automorphism group of
$M(1)$ is exactly $O(\h)$ (cf. \cite[Proposition 2.3]{DM}).

(b) Vertex operator algebra associated to the highest weight
modules for affine Kac-Moody algebras (cf.  \cite{DL}, \cite{FZ},
\cite{Li}). Let $\g=\h+\oplus_{\alpha\in \Delta}\g_{\alpha}$ be a
finite dimensional simple Lie algebra with a Cartan subalgebra
$\h$ and the corresponding root system $\Delta.$ We fix the
positive roots $\Delta_+$ and assume that $\theta$ is the longest
positive root. Denote $P_+$ the set of dominant weights. Let
$(\cdot,\cdot)$ be a nondegenerate symmetric invariant bilinear
form on $\g$ such that $(\theta,\theta)=2.$

The affine Kac-Moody algebra is
$$\hat{\g}={\g}\otimes\C[t,t^{-1}]\oplus \C K,$$
with bracket
$$[a(m),b(n)]=[a,b](m+n)+m\delta_{m+n,0}(a,b)K,$$
$$[K,\hat\g]=0,$$
where $a(m)=a\otimes t^m$ for $a\in \g$ and $m\in \Z.$ If $M$ is
an irreducible $\hat\g$-module then the center $K$ acts as a
constant $k$ which is called the {\em level} of the module.

Let $M$ be a $\g$-module and $k$ a complex number.
 The generalized Verma modules of level
 $k$ associated to $M$ is
$$\hat M_k=U(\hat\g)\otimes_{U(\g\otimes\C[t]\oplus \C K)}M$$
where $\g\otimes t\C[t] \cdot M=0,$ $\g\otimes t^0$ acts on $M$ as
$\g$ and $K=k$ on $M.$

Let $L(\lambda)$ be the  highest weight module for $\g$ with
highest weight $\lambda\in\h.$ Set
$V(k,\lambda)=\widehat{L(\lambda)}_k$ and denote $L(k,\lambda)$
the unique irreducible quotient. Then $L(k,\lambda)$ is integrable
or unitary if and only if
 $k$ is a nonnegative integer, $\lambda\in P_{+},$ and
$(\lambda,\theta)\leq k.$

Denote the dual Coxeter number of $\g$ by $h^{\vee}.$ Then
 $h^{\vee}$ can be defined by  $\sum_{\alpha\in\Delta}\frac{(\alpha,\alpha)}{2}=d h^{\vee}$
 where $d$ is the rank of $\g.$
Then we have

(a) Any $\hat\g$-quotient module of $V(k,0)$ is a vertex operator
algebra, if $k+h^{\vee}\ne 0.$

(b)$L(k,0)$ is rational if and only if $k$ is a nonnegative
integer. In this case, the irreducible $L(k,0)$-modules are
exactly the level $k$ unitary highest weight modules. In
particular, if $\g$ is the Lie algebra of type $D_d,$ then $\g$ is
the affine Lie algebra $D_{d}^{(1)}$ which has exactly 4 level 1
unitary highest weight modules. These modules are used in the
construction of elliptic genera.

The automorphism groups of $V(k,0)$ and $L(k,0)$  are exactly the
automorphism group of the Lie algebra $\g$ (cf. \cite{DM}).

\section{\normalsize  Vertex operator algebras bundles } \label{s3}
\setcounter{equation}{0}

$\quad$ In this section we use the category $\O_V$ to define
vertex operator algebra bundles over a compact space. This
generalizes the notion of complex vector bundles. The main result
in this section is that for any $V$-bundle $E$ there is another
$V$-bundle $F$ such that $E\oplus F$ is a trivial bundle. As in
the classical case, this result is crucial in defining VOA
$K$-groups and studying the cohomology properties of the
$K$-groups. We also prove that a $V$-bundle is essentially a
vector bundle if $V$ is rational.

\subsection{ \normalsize Definition of a VOA bundles and some consequences}
$\quad$ We now fix a vertex operator algebra $V$. Let $X$ be a
compact topological space.

\begin{defn}\label{t3.1}  Let $M\in \O_V.$ A $V$-bundle $E$ over $X$
with fiber $M$ is a direct sum of vector bundles $E=\oplus_{\l\in
\C}E^\l$ over $X$ such that all transition functions are
$V$-module isomorphisms. That is, there exists an open covering
$\{X_\a\}_{\a\in I}$ of $X$ and  a family of continuous
isomorphism of vector bundles
 $H_\a=(H_\a^\lambda : E^\l|_{X_\a} \to X_\a\times M_\l )_{\l\in \C}$
 with $M=\oplus_{\l\in \C} M_\l$ a $V$-module such that if we denote by
$(H_\a^\lambda\circ (H_\b^\lambda)^{-1})_{\l\in \C} =
(g_{\a\b}^\l)_{\l\in \C}$, then $g_{\a\b}(x)=
(g_{\a\b}^\l(x))_{\l\in \C}: M\to M $ is a $V$-module isomorphism
for any $x\in  X_\a\cap  X_\b$.
\end{defn}

\begin{rem}\label{rem3}{\rm  (1) If $V=\C$ as in Remarks \ref{rem1}, the $V$-bundle defined here is exactly the classical complex vector bundle
over $X.$

(2) Let $U$ and $V$ be two vertex operator algebras and $E,F$ be
the $U$ and $V$-bundles over $X$ respectively. Then $E\otimes F$
is a $U\otimes V$-bundle over $X$ (see Remark \ref{raa}).
 In particular, If $U=\C$ then $E\otimes F$ again is a $V$-bundle over $X.$

(3) One could use the full module category of $V$-modules to
define  a $V$-bundle but the property (2) in this remark would not
be true. But this property is fundamental if one wants to
construct a ring instead of a group from vertex operator algebra
bundles (see Corollary \ref{c}).   }

(4) One can give a different definition of $V$-bundle using a
subgroup $G$ of $\Aut(V).$ In this case, we denote $\O_{V,G}$ the
subcategory  of $\O_V$ consisting of $V$-modules $M$ such that $G$
acts on $M$ as automorphisms. The transition functions
$g_{\a,\b}(x)$ now are required to be in $G.$ Most results in this
paper hold for this definition of $V$-bundles.

(5) There was a notion of vertex operator algebra bundle given in
\cite{FBZ}. But our definition is different from theirs.
\end{rem}

Sub-bundles and quotient bundles, direct sum of bundles are
defined as expected. Let $E,F$ be two $V$-bundles on $X$, a map
$f: E\to F$ is called a $V$-bundle homomorphism if there exist a
family of continuous homomorphisms of vector bundles $f^\l: E^\l
\to F^\l$ such that $f= (f^\l)_{\l\in \C} : E\to F$ is a
$V$-module homomorphism for any $x\in X$. We call a $V$-bundle
homomorphism $f: E\to F$ an  epimorphism, monomorphism and
isomorphism of $V$-bundles if $f^\l$ is so for any $\l\in \C$.

We say a $V$-bundle $E$ is trivial if there exists a $V$-bundle
isomorphism
 $\varphi:E\to X\times M$, here $X\times M$ is the natural
$V$-bundle on $X$ with the $V$-module $M$ as fibers.

We now discuss some well known examples of vertex operator algebra
bundles which have been used in the literature to study elliptic
genus and Witten genus.

(a) If $X$ is a Riemannian manifold, then the transition functions
of the complex tangent bundle $T_{\C}X$ lie in the special
orthogonal group $SO(d)$ where $d$ is the dimension of $X.$ Then
$\bigotimes_{n> 0}S_{q^n}(T_\C X)$ is a $M(1)^{SO(d)}$-bundle.
Here and below, $S_t(E)= 1+t\,E +t^2 \sym^2(E)+\cdots,$
$\Lambda_t(E)=1+t\,E +t^2 \Lambda^2(E)+\cdots$ are respectively
the symmetric and wedge operation of a vector bundle $E$. $M(1)$
is the Heisenberg vertex operator algebra of dimension $d$ with
$SO(d)$ as a subgroup of $\Aut(M(1))$ and $M(1)^{SO(d)}$ is the
$SO(n)$-invariants of $M(1)$ which is a vertex operator subalgebra
of $M(1).$ Similarly, $\bigotimes_{n\geq
0}\Lambda_{q^{n+1/2}}(T_\C X)$ is an $L(1,0)^{SO(d)}$-bundle where
$L(1,0)$ is the level one module for the affine algebra
$D_{d/2}^{(1)}$. Here we assume that $d$ is even in this case.

(b) If $X$ is further assumed to be a spin manifold, we denote the
spin bundle by $S.$ Then $S\bigotimes
\bigotimes_{n>0}\Lambda_{q^{n}}(T_\C X)$ is also a
$L(1,0)^{SO(d)}$-bundle.

\bigskip

Let $E$ be a $V$-bundle over $X.$ Set $E'=\oplus_{\l\in
\C}(E^\l)^*.$ Then $E'$ is also a $V$-bundle in an obvious way.

\begin{lem}\label{l3.1} Let $E$ be a $V$-bundle on $X$. Then

(1) $E\oplus E'$ is a $V$-bundle with a nondegenerate symmetric
invariant bilinear form induced from the natural bilinear from on
$M\oplus M'.$ In particular, any $V$-bundle on $X$ can be embedded
into a bundle with a nondegenerate symmetric invariant bilinear
form.

(2) Let $\{g_{\a\b}^*\}$ be the transition functions of $E'.$ Then
we have
$$(g_{\a\b}^*(x)s^*,g_{\a\b}(x)s)=(s^*,s)$$
for any $\a,\b\in I,$ $x\in X_\a\cap X_\b$ and $s\in M,$ $s^*\in
M'.$
\end{lem}

\pf (1) follows from Lemma \ref{l2.1} and (2) is clear. \qed

\begin{prop}\label{p3.4} For any $V$-bundle $E$ there exists another
$V$-bundle $F$ such that $E\oplus F$ is a trivial $V$-bundle.
\end{prop}
\pf By Lemma \ref{l3.1} we can assume that there is a
nondegenerate invariant symmetric bilinear form on $E.$

Certainly we can take a finite covering $\{X_\a\}_{\a\in I}$ in
the Definition \ref{t3.1}. Write $H_\a=(\pi,h_{\a})$ for
$\alpha\in I$ where $\pi: E \to X$ is the natural projection map
and $h_\a(e)\in M$. Let $\{p_\a : X\to \R_+ \}_{\a\in I}$ be a
family of continuous functions on
 $X$ such that $\{ p_\a ^2\}_{\a\in I}$ is a partition of unity on $X$
associated to the covering $\{X_\a\}_{\a\in I}$.
  Define a $V$-bundle homomorphism
$$\sigma: E\to X\times M^{\oplus n}$$
by $\sigma(e)=(x, (p_\a(x)h_{\a}(e))_{\a\in I})$ where $x=\pi (e)$
and $n$ is the cardinality of $I.$ Since $p_\a(x)>0$ for some
$\alpha\in I,$ we see that $\sigma$ is a monomorphism. So we can
identify $E$ with its image $\sigma(E)$ in the trivial $V$-bundle
$X\times  M^{\oplus n}.$

By extending the bilinear form from $M$ to $M^{\oplus n}$
componentwisely, we have a nondegenerate invariant symmetric
bilinear form on $X\times  M^{\oplus n}.$ By Lemma \ref{l3.1}, the
transition functions preserve the bilinear form on $M$, thus for
any $e,f\in E_x$
\begin{eqnarray*}
& & (\sigma(e),\sigma(f))_x
=\sum_{\alpha\in I}\Big  (p_\a(x) h_\a(e),p_\a(x) h_\a(f)\Big )\\
& &\ \ \ \ \ =\sum_{\alpha\in I}\Big  (g_{\a\b}(x) h_\b(e),
g_{\a\b}(x)h_\b(f) \Big )p_\a ^2(x)\\
& &\ \ \ \ \ =(h_\b(e),h_\b(f))\sum_{\alpha\in I}p_\a ^2(x)\\
& &\ \ \ \ \ =(e,f)_x,
\end{eqnarray*}
where $\b\in I$ is fixed such that $x\in X_\b.$ Thus the
homomorphism $\sigma$ preserves the bilinear form and the
restriction of the  bilinear form to $\sigma(E)$ is nondegenerate.
Set $F=\sigma(E)^{\perp}.$ Then $F$ is still a $V$-bundle on $X$
and $X\times M^{\oplus n}=\sigma(E)\oplus F\cong E\oplus F,$ as
desired. \qed

\subsection{ \normalsize The structure of $V$-bundles for rational VOA $V$}

$\quad$ Next we discuss a special case when $V$ is a rational
vertex operator algebra. Then $V$ has only finitely many
irreducible modules up to isomorphisms, say $\{M^1,...,M^p\}$. Any
$V$-module $W$ is isomorphic to $\oplus_{i=1}^pn_iM^i.$ The
nonnegative integer  $n_i$ is called the index of $M^i$ in $W$ and
is denoted by $[W:M^i].$ Given another $V$-module $M$ we define
the Grassmannian variety $G(W,M)$ of $W$ to be the set of
$V$-submodules of $W$ isomorphic to $M.$ It is clear that if
$[M:M^i]> [W:M^i]$ for some $i,$ $G(W,M)$ is an empty set. So we
restrict ourselves to the case when $M$ is a submodule of $W.$ Let
$W=\bigoplus_{i=1}^pk_iM^i$ with $k_i\leq n_i.$ Then
$$G(W,M)\cong G(n_1,k_1)\times \cdots \times G(n_p,k_p)$$
as varieties where $G(n,k)$ is the classical Grassmann variety.
That is $G(n,k)$ is the set of $k$-dimensional subspaces of an
$n$-dimensional complex vector space. Clearly, $G(W,M)$ is compact
as each $G(n_i,k_i)$ is compact.

Note that over $G(W,M)$ there is a canonical $V$-bundle $E(W,M)$
whose total space consists of $(H,x)$ where $H$ is a $V$-submodule
of $W$ isomorphic to $M$ and $x\in H.$

\begin{prop}\label{p3.5} Let $E$ be a $V$-bundle over $X$ whose fiber is isomorphic
to a $V$-module $M.$ Then there is a bundle homomorphism $f$ from
$E$ to $E(M^{\oplus n},M)$ such that $E$ is isomorphic to
$f^*E(M^{\oplus n},M).$
\end{prop}

\pf Recall that $\{X_\a|\alpha\in I\}$ is an open covering of $X$
and $\{p_\a^2\}_{\a\in I}$ is a partition of the unity with $p_\a$
being supported by $X_\a.$ Also recall the isomorphism
$$H_\a: E|_{X_\a}\to X_\a\times M$$
which sends $e$ to $(\pi(e), h_\a(e)).$

In order to construct the bundle homomorphism $f$ it is enough to
define a $V$-module monomorphism $\bar f: E_x \to M^{\oplus n}.$
Let $e\in E_x.$ As in the classical case, we set $\bar
f(e)=(p_\a(e) h_\a(e))_{\a\in I}\in M^{\oplus n}.$ It is clear
$\bar f$ is a $V$-module monomorphism. As usual, $f(e)=(\bar
f(E_x), \bar f(e)).$ Then $f$  is a bundle homomorphism from $E$
to $E(M^{\oplus n},M)$ such that $E$ is isomorphic to
$f^*E(M^{\oplus n},M).$ \qed\\

The result in Proposition \ref{p3.5} is not surprising. Such a
result is well known for ordinary vector bundles. But if $V$ is
rational, a $V$-bundle is determined by certain vector bundles. In
fact let $E$ be a $V$-bundle with fiber $M.$ Let
$M=\oplus_{i=1}^pW_i\otimes M^i$ where as before $\{M^1,...,M^p\}$
is a complete list of inequivalent irreducible $V$-modules and
$W_i$ is the space of multiplicity of $M^i$ in $M.$ Then each
$g_{\a\b}$ defines a map $h_{\a\b}:  X_\a\cap  X_\b \to
\bigoplus_{i=1}^p\End (W_i).$ For each $i$ we define a vector
bundle $V(E)^i$ over $X$ with fiber $W_i$ and transition functions
$\{h_{\a\b}|\a,\b\in I\}.$ For any $V$-module $N\in \O_V$ we
denote the trivial $V$-bundle on $X$ by $\underline{N}$. That is,
$\underline{N}=X\times N.$ From Remark \ref{rem3}, we have a
$V$-bundle $V(E)^i\otimes \underline{M}^i.$ So we have proved the
following proposition.
\begin{prop}\label{p3.6} If $V$ is rational then for any
$V$-bundle $E$ over $X$ there are vector bundles $V(E)^i$ for
$i=1,...,p$ such that $E$ is isomorphic to
$\oplus_{i=1}^pV(E)^i\otimes \underline{M}^i.$
\end{prop}

\section{\normalsize  Definition of VOA $K$-group} \label{s4}
\setcounter{equation}{0}

$\quad$ In this section we define a $K$-group associated to
$V$-bundles for a vertex operator algebra $V$. We follow closely
the set-up of \cite{A} with certain necessary extensions.

\subsection{\normalsize  Definition of $K_V(X)$} \label{s4.1}

$\quad$ Let $X$ be a compact space, we denote by ${\cal V}_V(X)$
the set of isomorphic classes of $V$-bundles over $X.$ ${\cal
V}_V(X)$ is an abelian semigroup with addition given by the direct
sum. We also denote by $K_V(X)$ the abelian group generated by the
equivalence classes of $V$-bundle $[E].$ Then the elements of
$K_V(X)$ are of the form $[E]-[F].$

\begin{rem}\label{rem4}{\rm If $V=\C,$ the $K_V(X)$ is the classical
$K(X)$ as defined in \cite{A}.}
\end{rem}

First from Remark \ref{rem3} we have the following:
\begin{lem}\label{l4.0} If $U,V$ are two vertex operator algebras then the tensor product
of the bundles induces a natural group homomorphism
$K_U(X)\otimes_{\Z} K_V(X)\to K_{U\otimes V}(X).$ In particular,
if $U=\mathbb{C},$ then $K_V(X)$ is a natural $K(X)$-module (cf.
Remark \ref{rem4}).
\end{lem}

The following corollary is an immediate consequence of Lemma
\ref{l4.0}.
\begin{cor}\label{c} For any vertex operator algebra $V$,
$$\oplus_{n\geq 0}K_{V^{\otimes n}}(X)$$
is a commutative algebra over $K(X)$ where $V^{\otimes 0}$ is
understood to be $\C.$
\end{cor}

The following proposition is crucial to prove the cohomological
properties of $K_V(X)$. We follow the argument of
\cite[Proposition 1.7]{H} which  avoid the use of the extension of
the section of $\Aut_V (E)$ from a compact subset to a open
neighborhood, as $\Aut_V (E)$ may be infinite dimensional.

\begin{lem} \label{l4.1} Let $Y$ be a compact Hausdorff space, $f_t: Y\to X$
$(0\leq t \leq 1)$ a homotopy and $E$ a $V$-bundle over $X$. Then
\begin{align}
f^*_0 E \simeq f^*_1 E.
\end{align}
\end{lem}
\pf  Denote by $I$ the unit interval and let $f: Y\times I \to X$
be the homotopy, so that $f(y,t)=f_t(y)$, and let $\pi: Y\times I
\to Y$ denote the projection onto the first factor.

  At first, we can choose a finite  open covering $\{Y_{x_i}\}_{i=1}^n$
of $Y$ so that
 $f^*E$ is trivial over each $Y_{x_i} \times I$. In fact, for each $x\in Y$
we can find open neighborhood $U_{x,1}, \cdots, U_{x,k}$ in $Y$
and a partition $0=t_0< t_1 \cdots < t_k=1$ of $[0,1]$ such that
the bundle is trivial over each $ U_{x,i}  \times
[t_{i-1},t_{i}]$. Set $Y_x = U_{x,1}\cap  \cdots \cap U_{x,k}$.
Now we claim that $f^*E$ is trivial over  $Y_x \times I$. To see
this, let $h_i: f^*E|_{Y_x \times [t_{i-1},t_{i}]} \to Y_x \times
[t_{i-1},t_{i}] \times M$ be the trivializing isomorphisms. We
take the $V$-bundle isomorphism $h'_1 (y,t) = (h_0\circ
h_1^{-1})(y,t_1) \circ h_1 (y,t): f^*E|_{Y_x \times
[t_{1},t_{2}]}\to Y_x \times [t_{1},t_{2}] \times M$, then
$h_1=h'_2$ on $Y_x \times \{t_1\}$, thus they define a
trivialization on $Y_x \times [t_{0},t_{2}]$, and in this way, we
know that $f^*E$ is trivial over  $Y_x \times I$. Now, as $Y$ is
compact, there exist $\{Y_{x_i}\}_{i=1}^n$ which cover $Y$.

Let $p_i$ be a partition of unity of $Y$ with support of $p_i$
contained in $Y_{x_i}$. For $i\geq 0$, set $q_j= \sum_{i=1}^j
p_i$. In particular $q_0=0$ and $q_n=1$. Let $W_i$ be the graph of
$q_i$, the subspace of $Y\times I$ consisting of points of the
form $(x,q_i(x))$, and let $\pi_i: E_i \to W_i$ be the restriction
of the bundle $E$ over $W_i$. Since $E$ is trivial on
$Y_{x_i}\times I$, the natural projection homeomorphism $W_i\to
W_{i-1}$ lifts to a homeomorphism $g_i : E_i\to E_{i-1}$ which is
identity outside  $\pi_i (Y_{x_i})$ and which takes each fiber of
$E_i$ isomorphically onto the corresponding fiber of $E_{i-1}$.
The composition $g= g_1\circ g_2\cdots \circ g_n$ is then an
isomorphism from the restriction of $E$ over $Y\times \{1\}$ to
the restriction on $Y\times \{0\}$. \qed

\begin{lem}\label{l4.3}  (1) Any element of  $K_V(X)$ can be represented by an element
of the form $[E]-[\underline{M}],$ where $E$ is a $V$-bundle and
$M$ is a $V$-module.

(2) If $[E]=[F]$ in $K_V(X)$ then there is a $V$-module $M$ such
that $E\oplus \underline{M}\cong F\oplus\underline{M}.$
\end{lem}

\pf As we have mentioned already, every element of $K_V(X)$ is of
the form $[H]-[G].$ By Proposition \ref{p3.4}, there exists a
$V$-bundle $F$ and a $V$-module $M$ such that $H\oplus F\cong \M.$
Thus we have
$$[G]-[H]=[G+F]-[H+F]=[G+F]-[\M].$$
This proves (1).

If $[E]=[F],$ then there exists a $V$-bundle $G$ such that
$E\oplus G\cong F\oplus G.$ Let $H$ be a $V$-bundle such that
$G\oplus H\cong \M$ for some $V$-module $M.$ Then we have $E\oplus
\M\cong F  \oplus \M.$ This proves (2). \qed

$\ $

 Clearly, $K_V(pt)$ is just
the the Grothendieck group of  $V$-modules. From Proposition
\ref{p3.6} we have the following proposition.

\begin{prop}\label{l4.4} If $V$ is rational, then
\begin{align}
K_V(X)=K(X)\otimes_{\Z} K_V(pt).
\end{align}
\end{prop}

\begin{rem} {\rm If $V$ is
rational and  $\{M^1,...,M^p\}$ is the set of irreducible
inequivalent $V$-modules up to isomorphisms, then $K_V(pt)$ is
isomorphic to the group $\Z^p$ with generators $[M^1],...,[M^p].$
In particular, for $V$-modules $M$ and $N,$ $[M]=[N]$ if and only
if $M\cong N$ as $V$-modules.}
\end{rem}

\subsection{\normalsize  Definition of $K_V(X,Y)$}\label{s4.2}

$\quad$ We next define $K_V(X,Y)$ for a compact pair $(X,Y)$. Let
$\cC$ denote the category of compact spaces, $\cC^+$ the category
of compact spaces with distinguished basepoint, and $\cC^2$ the
category of compact pairs. We define a functor $\cC^2 \to \cC^+$,
 by sending a pair $(X,Y)$
to $X/Y$ with base point $Y/Y$ (if $Y\neq \emptyset$, the empty
set, $X/Y$ is understood to be the disjoint union of $X$ with a
point.). If $X\in\cC$, we denote $(X,\emptyset)\in \cC^2$  by
$X^+$.

If $X$ is in $\cC^+$, we define $\wi{K}_V(X)$ to be the kernel of
 the map $i^*: K_V(X) \to K_V(x_0)$ where $i:x_0\to X$ is the
inclusion of the basepoint. If $c: X\to x_0$ is the collapsing map
then $c^*$ induces a splitting $K_V(X)=\wi{K}_V(X)\oplus
K_V(x_0)$.
 This splitting is clearly natural for maps in $\cC^+$.
Thus $\wi{K}_V$ is a functor on $\cC^+$. Also, it is clear that
$K_V(X)\simeq \wi{K}_V(X ^+)$. We define $K_V(X,Y)$ by $K_V(X,Y)=
\wi{K}_V(X/Y)$. In particular $K(X,\emptyset) \simeq K_V(X)$.
Since $\wi{K}_V$ is a functor on $\cC^+$ it follows that
$K_V(X,Y)$ is a contravariant functor of $(X,Y)$ in $\cC^2$.

We now introduce the smash product operator in $\cC^+$, if $X,Y\in
\cC^+$, we put $X\wedge Y = X\times Y /X \vee Y $ where $X \vee Y
= X\times \{y_0\} \cup \{x_0\} \times Y$, $x_0, y_0$
 being the base-points of $X,Y$ respectively. For any three spaces
$X,Y,Z\in \cC^+$, we have a natural homeomorphism $X\wedge
(Y\wedge Z)\simeq (X\wedge Y)\wedge Z$
 and we shall identify these spaces by the homeomorphism.

Let $I$ denote the unit interval $[0,1]$ and let $\partial I =
\{0\}\cup \{1\}$ be its boundary. We take $I/\partial I\in \cC^+$
as our standard model of the circle $S^1$. For $X\in \cC^+$ the
space $S^1\wedge X\in \cC^+$ is called the reduced suspension of
$X$, and often written as $SX$. The $n$-th iterated suspension
$SS\cdots SX$ ($n$ times) is naturally homeomorphic to $S^n \wedge
X$ and is written briefly as $S^nX$.

\begin{defn}\label{t4.5}  For $n\leq 0$,
\begin{align}\label{}
& \wi{K}^{-n}_V (X) = \wi{K}_V(S^n X)\ {\rm for}\  X\in \cC^+,\\
&K^{-n}_V(X,Y)=  \wi{K}^{-n}_V (X/Y)=
\wi{K}_V(S^n (X/Y))\  {\rm for}\  (X,Y)\in \cC^2,\nonumber  \\
&{K}^{-n}_V (X) = K^{-n}_V(X,\emptyset) = \wi{K}_V(S^n (X^+)) \
{\rm for}\  X\in \cC.\nonumber
\end{align}
\end{defn}

By proceeding further, we define the cone on $X$ by
$$ CX= I\times X/\{0\} \times X.$$
Thus $C$ is a functor $C: \cC \to \cC^+$. We identify $X$ with the
subspace $\{1\}\times X$ of $CX$. The space $CX/X = I\times X
/\partial I\times X$ is called the {\em unreduced suspension} of
$X$.

\section{\normalsize  Cohomological properties of $K$}\label{s5}
\setcounter{equation}{0}

$\quad$ In this section, we assume that $X$ is a finite CW-complex
 and  $Y\subset X$ is a CW sub-complex. This condition allows us to extend a
 trivialization of $V$-bundles on $Y$ to a neighborhood of $Y$ which
 is crucial in the proof of Lemma \ref{t5.1}, Theorem \ref{t5.3}.
We will establish the basic
 cohomological properties of $K_V$ for these compact pairs $(X,Y)$.

\begin{lemma}  \label{t5.1} We have an exact sequence
\begin{eqnarray} \label{5.1}
K_V(X,Y) \stackrel{j^*}{\to} K_V(X) \stackrel{i ^*}{\to} K_V(Y),
\end{eqnarray}
where $i: Y\to X$ and $j: (X,\emptyset)\to (X,Y)$ are the
inclusions.
\end{lemma}

\pf  The composition $i ^*\circ j^*$ is induced by the composition
$j\circ i: (Y,\emptyset)\to (X,Y)$ and so factor through the zero
group. Thus $i^*\circ j^*=0$. Suppose now that $\xi\in {\rm Ker}
i^*$. We may represent $\xi$ in the form $[E]-[\M]$ where $E$ is a
$V$-bundle over $X$ and $M$ is a $V$-module. Since $i ^* \xi =0$,
it follows that $[E]|_{Y} = [\M]$ in $K_V(Y)$. This implies that
there exists a $V$-module $N$ such that we have
\begin{eqnarray} \label{5.2}
 (E\oplus \underline{N})|_Y = \M\oplus \underline{N}.
\end{eqnarray}
Now as $Y$ is a CW sub-complex of $X$, there exists an open
neighborhood $U$ of $Y$ in $X$ such that $Y$ is a strong
deformation retract of $U$, i.e. there exists $f_t: U\to U$ ($t\in
[0,1]$) such that $f_1={\rm Id}_U$, $f_0|_Y={\rm Id}_Y$ and
$f_0(U)=Y$. By Lemma \ref{l4.1}, (\ref{5.2}), $(E\oplus \underline{N})|_U$ is
trivial on $U$. This defines a bundle $E\oplus \underline{N}/\alpha$ on $X/Y$
and so an element
\begin{eqnarray}
\tau= [E\oplus \underline{N}/\alpha] - [\M\oplus \underline{N} ] \in \wi{K}_V(X/Y) =
K_V(X,Y).
\end{eqnarray}
Then
\begin{eqnarray}
j^* (\tau) =[E\oplus \underline{N}/\alpha] - [\M\oplus \underline{N}] 
=  [E]-[\M]=\xi.
\end{eqnarray}
Thus ${\rm Ker} i ^* = {\rm Im} j ^*$ and the exactness is
established. \qed

\begin{cor} \label{t5.2} If $(X,Y) \in \cC ^2$ and $ Y\in \cC ^+$, then the
following sequence is exact,
\begin{eqnarray}
K_V(X,Y)\to \wi{K}_V(X) \to \wi{K}_V(Y).
\end{eqnarray}
\end{cor}
\pf This is immediate from Lemma \ref{t5.1} and the natural
isomorphisms
\begin{align}
 K_V(X)\simeq \wi{K}_V(X)\oplus K_V(y_0),\\
 K_V(Y)\simeq \wi{K}_V(Y)\oplus K_V(y_0). \nonumber
\end{align}
\qed

Our main proposition of this section is following,
\begin{thm}\label{t5.3} There is a natural exact
sequence (infinite to the left )
\begin{multline}\label{sequence0}
\cdots \to K^{-2}_V(Y)\stackrel{\delta}{\to}K^{-1}_V(X,Y)
\stackrel{j^*}{\to}K^{-1}_V(X) \stackrel{i ^*}{\to}  K^{-1}_V(Y)\\
\stackrel{\delta}{\to}K^{0}_V(X,Y)\stackrel{j^*}{\to}K^{0}_V(X)
\stackrel{i ^*}{\to}  K^{0}_V(Y).
\end{multline}
\end{thm}

\pf First we observe that it is sufficient to show that
\begin{align} \label{sequence1}
K^{-1}_V(X,Y) \stackrel{j^*}{\to}K^{-1}_V(X) \stackrel{i ^*}{\to}
K^{-1}_V(Y)
\stackrel{\delta}{\to}\wi{K}^{0}_V(X,Y)\stackrel{j^*}{\to}
\wi{K}^{0}_V(X) \stackrel{i ^*}{\to}  \wi{K}^{0}_V(Y)
\end{align}
is exact. In fact, if this has been established then, by replacing
$(X,Y)$ by $(S^nX, S^n Y)$ for $n=1,2,\cdots$, we obtain an
infinite sequence continuing (\ref{sequence1}). Then by replacing
$(X,Y)$ by $(X^+, Y^+)$ where $(X,Y)$ is any pair in $\cC^2$ we
get the infinite sequence of the enunciation. Now, Corollary
\ref{t5.2} gives the exactness of the last three terms of
(\ref{sequence1}). To get exactness at the remaining places we
shall apply  Corollary \ref{t5.2} in turn to the pairs $(X\cup
CY,X)$ and $ ((X\cup CY)\cup CX, X\cup CY)$. First, by taking the
pair $(X\cup CY, X)$, we get an exact sequence
\begin{align}\label{sequence2}
K_V(X\cup CY,X) \stackrel{\tau_1 ^*}{\to}  \wi{K}_V(X\cup CY)
\stackrel{\tau_2^*}{\to}  \wi{K}_V(X),
\end{align}
where $\tau_1, \tau_2$ are the natural inclusions. Let $U$ be the
neighborhood of $Y$ in $X$ as in Lemma \ref{t5.1}. Since $CY$ is
contractible, by Lemma \ref{l4.1}, any $V$-bundle $E$
 on $X\cup CY$ is trivial on $U\cup CY$, thus
$p^*:  \wi{K}_V(X/Y) \to \wi{K}_V(X\cup CY)$ is an isomorphism
where $p: X\cup CY \to X\cup CY/CY =X/Y$ is the collapsing map.
Also the composition $\tau_2 ^*\circ p^*$ coincides with $j^*$.
Let $\theta: K_V(X\cup CY,X)\to K^{-1}_V(Y)$ be the isomorphism
introduced earlier. Then by defining $\delta: K^{-1}_V(Y) \to
K_V(X,Y)$
 by $\delta = \tau_1^* \circ \theta ^{-1}$, we obtain the exact sequence
\begin{align}\label{sequence3}
K^{-1}_V(Y)\stackrel{\delta}{\to} K_V(X,Y) \stackrel{j
^*}{\to}\wi{K}_V(X),
\end{align}
which is the middle part of (\ref{sequence1}).

Finally, we apply  Corollary \ref{t5.2} to the pair $(X\cup C_1
Y\cup C_2 X, X\cup C_1 Y)$, where we have labeled the cones $C_1$
and $C_2$ in order to distinguish between them. Thus we obtain an
exact sequence
\begin{align}\label{sequence4}
K_V(X\cup C_1 Y\cup C_2 X, X\cup C_1 Y)\to \wi{K}_V(X\cup C_1
Y\cup C_2 X) \to \wi{K}_V(X\cup C_1 Y).
\end{align}
It will be sufficient to show that this sequence is isomorphic to
 the sequence obtained from  the first three terms  of (\ref{sequence1}).
In view of the definition of $\delta$, it will be sufficient to
show
 that the following diagram  commutes up to sign.
\begin{align}\label{sequence5}
 &K_V(X\cup C_1 Y\cup C_2 X, X\cup C_1 Y)& = \wi{K}_V(C_2X/X)
= &K^{-1}_V(X)\\
 &\stackrel{\downarrow}{\wi{K}_V(X\cup C_1 Y\cup C_2 X)}& =\wi{K}_V(C_1Y/Y) =& \stackrel{i^* \downarrow}{K^{-1}_V(Y)}.\nonumber
\end{align}

As in \cite[p74]{A}, we will get (\ref{sequence5}) if we can prove
that the following diagram commutes up to sign,
\begin{align}\label{sequence6}
K_V(C_1 Y\cup C_2 Y) \leftarrow &K_V(C_1 Y/Y) \leftarrow   \wi{K}_V(SY)\\
\stackrel{\nwarrow}{ } &K_V(C_2 Y/Y) \leftarrow
\stackrel{\parallel}{\wi{K}_V(SY)} . \nonumber
\end{align}
This follow from the following lemma \ref{t5.4}. \qed

\begin{lemma}\label{t5.4} Let $T: S^1\to S^1$ be defined by $T(t)= 1-t$,
$t\in I=[0,1]$.  Recall that $S^1= I/\partial I$. Let $T\wedge 1 :
SY \to SY$ be the map induced by $T$ on $S^1$ and the identity on
$Y$ for $Y\in \cC^+$. Then $(T\wedge 1)^* a = -a$ for $a\in
\wi{K}_V(SY)$.
\end{lemma}

\pf  By the construction and Lemma \ref{l4.3}, for any $a\in
\wi{K}_V(SY)$, there exist a $V$-module $M$ and a $V$-bundle $E$
on $SY$ such that $a= [E]- [\M]$. We define $E$ first. As $SY= C_1
Y\cup_Y C_2 Y$, and $ C_i Y (i=1,2)$ is contractible, we know that
$E|_{ C_i Y}$ are trivial, thus there are maps $f_i :E|_{ C_i Y}
\to  C_i Y \times M$ of $V$-module isomorphism. The composite $f=
f_2 \circ f_1^{-1}=(f^\l)_{\l\in \C}: Y \to {\Aut}_V (M)$ is well
defined. Now, the operation $(T\wedge 1)^*$ on  $[E]$ corresponds
to the operation   of replacing the map $y\to f(y)$ by $y\to
f(y)^{-1}=((f^\l)^{-1})_{\l\in \C}.$ We denote the corresponding
bundle by $E_1 \in {K}_V(SY)$. We need to prove that in
${K}_V(SY)$,
\begin{align}\label{sequence7}
[E]\oplus [E_1]= [\M] \oplus [\M].
\end{align}

 For $0\leq t\leq \pi/2$, $y\in Y$, set
\begin{align}\label{sequenc8}
F_t(y) = \begin{pmatrix} f(y)& 0 \\ 0&  1 \end{pmatrix}
 \begin{pmatrix}\cos (t) & \sin (t) \\ -\sin (t)& \cos (t) \end{pmatrix}
\begin{pmatrix} 1 & 0 \\ 0& f(y)^{-1} \end{pmatrix}
\begin{pmatrix} \cos (t) & -\sin (t) \\ \sin (t)& \cos (t) \end{pmatrix}.
\end{align}
Then
\begin{align}\label{sequenc9}
F_{\pi/2}(y)= \begin{pmatrix} 1&0 \\ 0&1 \end{pmatrix},\quad
F_0(y) = \begin{pmatrix} f(y)&0  \\ 0& f(y)^{-1} \end{pmatrix}.
\end{align}
This means that $F_t$ is a homotopy from $\begin{pmatrix} 1&0 \\
0&1 \end{pmatrix}$ to $ \begin{pmatrix} f(y)&0  \\0&  f(y)^{-1}
\end{pmatrix}$. Thus we get (\ref{sequence7}) from Lemma
\ref{l3.1}. \qed

Now, by Lemma \ref{l4.1} and Theorem \ref{t5.3}, we get

\begin{cor}\label{t5.5} If $Y$ is a retract of $X$, then for all $n\leq 0$, the sequence
$K^{-n}_V (X,Y) \to K^{-n}_V (X) \to K^{-n}_V (Y)$ is a split
short exact sequence, and
\begin{align}\label{sequenc10}
K^{-n}_V (X)\simeq K^{-n}_V (X,Y)\oplus K^{-n}_V (Y).
\end{align}
\end{cor}

\section{Associative algebra bundles}

$\quad$ This section is motivated by the relation between a vertex
operator algebra $V$ and its Zhu's algebra $A(V).$ It turns out we
can define associative algebra bundles for a large class of
associative algebras.

\subsection{\normalsize  Definition of associative algebra bundles} \label{s6.1}

$\quad$ We assume that $A$ is an associative algebra over $\C$
with an anti-involution $\sigma.$ The setting and results in
Sections 2-5 hold in the present situation with suitable
modifications.

Let $M$ be an $A$-module and we denote the dual space of $M$ by
$M'$ as before. We also denote the natural pairing $M'\times M\to
\C$ by $(m',m)$ for $m'\in M'$ and $m\in M.$  The following lemma
is obvious.

\begin{lem}\label{l6.1} $M'$ is also a $A$-module such that
$(am',m)=(m',\sigma(a)m)$ for $a\in A,$ $m'\in M'$ and $m\in M.$
\end{lem}

As in Section 2 a form $(\ ,\ )$ on an $A$-module $W$ is called
invariant if $(aw_1,w_2)=(w_1,\sigma(a)w_2)$ for $w_i\in W$ and
$a\in A.$

\begin{lem}\label{l6.2} Lemma \ref{l2.1} holds for an $A$-module
$M.$
\end{lem}

We also need to define the category ${\O}_A$ of $A$-modules. An
$A$-module $W$ is in ${\O}_A$ if there exist $\l_1,...,\l_s\in\C$
such that $W=\bigoplus_{i=1}^s\bigoplus_{n\geq 0}W_{\l_i+n}$ is a
direct sum of finite dimensional $A$-modules and such that
$\Hom_{A}(W_{\l},W_{\mu})=0$ if  $\mu\ne \l.$  Such definition of
category $\O_A$ is well justified by Remark \ref{rem6.3}, Theorems
\ref{P3.1} (7) and \ref{t6} below.

\begin{defn}\label{d6.1}  Let $X$ be a compact space and
$W\in \O_A.$ An $A$-bundle $E$ over $X$ with fiber $W$ is a direct
sum of vector bundles $E=\oplus_{\l\in \C}E^\l$ over $X$ such that
all transition functions are $A$-module isomorphisms. That is,
there exists an open covering $\{X_\a\}_{\a\in I}$ of $X$ and a
family of isomorphism of vector bundles
 $H_\a=(H_\a^\lambda : E^\l|_{X_\a} \to X_\a\times M_\l )_{\l\in \C}$
 with $M=\oplus_{\l\in \C} M_\l$ an $A$-module such that if we denote by
$(H_\a^\lambda\circ (H_\b^\lambda)^{-1})_{\l\in \C} =
(g_{\a\b}^\l)_{\l\in \C}$,
 then each $g_{\a\b}^\l(x): M_\l\to M_\l$ is an
$A$-module isomorphism for any $x\in  X_\a\cap  X_\b,$ $\l\in\C.$
In particular, each $E^\l$ is an $A$-bundle.
\end{defn}

The analogues of Remark \ref{rem3} are as follows:
\begin{rem}\label{rem6.3}{\rm  (1) If $A=\C,$ then $\sigma$ is necessarily
the identity map,  the $A$-bundle defined here is exactly the
classical complex vector bundle over $X.$

(2) Let $A$ and $B$ be two associative algebras with
anti-involutions $\sigma_A$ and $\sigma_B$ respectively. Then
$A\otimes_\C B$ is an associative algebra with anti-involution
$\sigma_A\otimes \sigma_B.$
 Assume $E,F$ are
the $A$ and $B$-bundles over $X.$ then $E\otimes F$ is a
$A\otimes_{\C} B$-bundle over $X.$ In particular, If $A=\C$ then
$E\otimes F$ again is a $B$-bundle over $X.$}
\end{rem}

We can also define subbundles, quotient bundles, direct sum of
bundles. Various bundle homomorphisms are also defined as
expected.

Let $E$ be an $A$-bundle over $X.$ Regarding $E$ as a vector
bundle over $X,$ the dual bundle $E'$ is also an $A$-bundle.

Lemma \ref{l3.1} and  Proposition \ref{p3.4} hold with $V$
replaced by $A.$

We also define the K-group $K_A(X)$ to be the abelian group
generated by the equivalence classes of $A$-bundles.  Then Lemmas
\ref{l4.0}, \ref{l4.1}, \ref{l4.3} and Corollary \ref{c} hold with
obvious changes. Proposition \ref{l4.4} is also true if $A$ is
semisimple (that is, $A$ is a direct sum of full matrix algebras)
and with $V$ replaced by $A.$ That is $K_A(X)=K(X)\otimes_\Z
K_A(pt).$

We can also define $K_A(X,Y)$ and related objects as in Section 4
and the cohomological properties in Section 5 also hold in the
present setting.

\subsection{\normalsize  Zhu's algebra $A(V)$} \label{s6.2}

$\quad$ We review the Zhu's algebra $A(V)$ in this subsection
following \cite{Z} and \cite{DLM2}.

Let $V$ be a vertex operator algebra.  For any homogeneous vectors
$a\in V$, and $b\in V$, we define
\begin{align*}
&a*b=\left(\Res_{z}\frac{(1+z)^{\wt{a}}}{z}
Y(a,z)\right)b,\\
&a\circ b=\left( \Res_{z}\frac{(1+z)^{\wt{a}}}{z^{2}}
Y(a,z)\right)b,
\end{align*}
and extend to $V\times V$ bilinearly, here $\Res_z$ denotes the
coefficient of $z^{-1}$.  Denote by $O(V)$ the linear span of
$a\circ b$ ($a,b\in V$) and set $A(V)=V/O(V)$.  We write $[a]$ for
the image $a+O(V)$ of $a\in V$.

For homogeneous $a\in V$ we set $o(a)=a_{\wt a-1}$ and extend
linearly to all $a\in V.$ If $M=\bigoplus_{n\geq 0}M(n)$ is an
admissible module then $o(a)M(n)\subset M(n)$ for all $a\in V$ and
$n\in\Z.$ We now define the space of lowest weight vectors of $M:$
$$\Omega(M)=\{w\in M|a_{\wt a+m}w=0, a\in V, m\geq 0\}.$$

\begin{rem}{\rm If $M=\bigoplus_{\l\in\C}M_{\l}$ is a $V$-module
then $\Omega(M)=\bigoplus_{\l\in\C}\Omega(M)_{\l}$ is naturally
graded and each homogeneous subspace $\Omega(M)_{\l}=\Omega(M)\cap
M_{\l}$ is finite dimensional. }
\end{rem}

The following lemma is evident.
\begin{lem}\label{l6.a} If $M$ and $W$ are admissible $V$-modules and $f: M\to W$
is a $V$-module homomorphism then $f(\Omega(M))\subset \Omega(W).$
In particular, if $f$ is an isomorphism then
$f(\Omega(M))=\Omega(W).$
\end{lem}

The following theorem is due to \cite[\S 2]{Z} (also see
\cite{DLM2}).
\begin{thm}\label{P3.1}
(1) The bilinear operation $*$ induces on $A(V)$ an associative
algebra structure. The vector $[\1]$ is the identity and $[\w]$ is
in the center of $A(V)$.

(2) The linear map
$$\phi:  a\mapsto e^{L(1)}(-1)^{L(0)}a$$
induces an anti-involution of $A(V).$

(3) Let $M=\bigoplus_{n=0}^{\infty}M(n)$ be an admissible
$V$-module with $M(0)\ne 0.$ Then the linear map
\[
o:V\rightarrow\End (\Omega(M)),\;a\mapsto o(a)|\Omega(M)
\]
induces an algebra homomorphism from $A(V)$ to $\End (\Omega(M))$.
Thus $\Omega(M)$ is a left $A(V)$-module.

(4) The map $M\mapsto M(0)$ induces a bijection from the set of
equivalence classes of irreducible admissible $V$-modules to the
set of equivalence classes of irreducible $A(V)$-modules.

(5) If $M=\bigoplus_{\l\in\C}M_{\l}$ is a $V$-module, then each
$\Omega(M)_{\l}$ is an finite-dimensional $A(V)$-module.

(6) If $V$ is rational then $A(V)$ is a finite dimensional
semisimple algebra. Moreover, $M\to \Omega(M)$ gives an
equivalence of the category of admissible $V$-modules and the
category of $A(V)$-modules and the same functor also gives an
equivalence of the category of $V$-modules and the category of
finite dimensional $A(V)$-modules.

(7) If $M\in \O_V$ then $\Omega(M)\in \O_{A(V)}.$
\end{thm}

Note that if $M$ is a $V$-module and $\l\ne \mu$ then
$\Hom_{A(V)}(\Omega(M)_\l,\Omega(M)_{\mu})=0$ as $[\w]$ acts on
$\Omega(M)_{\mu}$ as scalar $\mu.$ Thus $\Omega(M)$ is an element
of $\O_{A(V)}.$ In fact, the definition of $\O_{A}$ for an
associative algebra $A$ reflects well the properties of
$\Omega(M).$

\subsection{\normalsize Relationship between $K_V(X)$ and $K_{A(V)}(X)$} \label{s6.3}

$\quad$ The relationship between $K_V(X)$ and $K_{A(V)}(X)$  is a
reflection of the relation between $V$ and $A(V)$ for any vertex
operator algebra $V$ and compact space $X.$

Recall Definition \ref{t3.1}. Let $E$ be a $V$-bundle over $X$
with fiber $M=\bigoplus_{\l\in \C}M_{\l}.$ Let $\{X_{\a}|\a\in
I\}$ be an open covering of $X$ which gives a local trivialization
of $E.$ We define a graded vector bundle
$\Omega(E)=\bigoplus_{\l\in\C}\Omega(E)^{\l}$ in the following
way: set
$\Omega(E)^{\l}|_{X_{\alpha}}=(H_{\alpha}^{\l})^{-1}(X_\a\times
\Omega(M)_\l)$ and $\Omega(H_{\alpha})^\l=H_{\a}|_{\Omega(E)^\l}.$
 Then $\Omega(H_{\alpha})^{\l}: \Omega(E)^{\l}|_{X_{\alpha}}\to X_{\alpha}\times \Omega(M)_{\l}$
  is an isomorphism. Set $\Omega(H_{\a})=(\Omega(H_\a)^\l)_{\l\in\C}$
and $\Omega(g_{\a\b})^{\l}=\Omega(H_\a)^\lambda\circ
(\Omega(H_\b)^\lambda)^{-1}$ for $\l\in \C.$

Then  $\Omega(g_{\a\b})(x)= (\Omega(g_{\a\b})^\l(x))_{\l\in \C}:
\Omega(M)\to \Omega(M)$ is an $A(V)$-module isomorphism  for any
$x\in  X_\a\cap  X_\b$ by Lemma \ref{l6.a}. By Theorem \ref{P3.1}
(7),  $\Omega(E)$ is an $A(V)$-bundle over $X.$

\begin{thm}\label{t6} The map $\Omega: E\to \Omega(E)$ from the set of $V$-bundles
over $X$ to the set of $A(V)$-bundles over $X$ induces a
homomorphism, which we still denote by $\Omega$, from the group
$K_V(X)$ to the group $K_{A(V)}(X).$ In particular, if $V$ is
rational, $\Omega$ is an isomorphism.
\end{thm}

\pf It is clear that $\Omega$ is a group homomorphism. Now we
assume that $V$ is rational. Recall Proposition \ref{l4.4} and its
analogue $K_{A(V)}(X)=K(X)\otimes_{\Z}K_{A(V)}(pt)$ for $A(V).$
Using Theorem \ref{P3.1} (5) we need to prove that $K_V(pt)$ and
$K_{A(V)}(pt)$ are isomorphic groups under the map $\Omega.$   But
this is clear from Theorem \ref{P3.1} (5) again by noting that the
category $\O_{A(V)}$ is exactly the category of finite dimensional
$A(V)$-modules. \qed

\end{document}